\documentclass[11pt]{amsart}
\usepackage{amsmath}
\usepackage{amsthm}
\usepackage{latexsym}
\usepackage{amssymb}
\usepackage{xspace}
\usepackage{amscd}

\newcommand{\Gal}{{\rm Gal}}
\newcommand{\lcm}{{\rm lcm}}

\newcommand{\inv}{^{-1}}

\newcommand{\pr}{{\mathbb P}}
\newcommand{\N}{{\mathbb N}}
\newcommand{\Z}{{\mathbb Z}}
\newcommand{\Q}{{\mathbb Q}}
\newcommand{\R}{{\mathbb R}}

\newcommand{\HQ}{{\mathbb H}}

\newcommand{\matriz}[1]{\begin{array} #1 \end{array}}

\newcommand{\GEN}[1]{\langle #1 \rangle}

\newcommand{\quat}[2]{\left( \frac{#1}{#2} \right)}

\title[The subgroup generated by cyclic cyclotomic algebras]
{The gap between the Schur group and the subgroup generated by cyclic cyclotomic algebras}

\author[A. Herman]{Allen Herman}
\address{Department of Mathematics, University of Regina, Canada}
\email{aherman@math.uregina.ca}
\author[G. Olteanu]{Gabriela Olteanu}
\address{Department of Mathematics and Computer Science, North University of Baia Mare,
Victoriei 76, 430122 Baia Mare, Romania. Current address. Departamento de Matem\'{a}ticas,
Universidad de Murcia, Espa\~{n}a}
\email{olteanu@math.ubbcluj.ro, golteanu@um.es}
\author[\'{A}. del R\'io]{\'{A}ngel del R\'io}
\address{Departamento de Matem\'{a}ticas, Universidad de Murcia, Espa\~{n}a}
\email{adelrio@um.es}

\thanks{{Research supported by the National Science and Engineering Research
Council of Canada, D.G.I. of Spain and Fundaci\'{o}n S\'{e}neca of Murcia.}}

\newtheorem{theorem}{Theorem}[section]
\newtheorem{lemma}[theorem]{Lemma}

\newtheorem{corollary}[theorem]{Corollary}
\newtheorem{remark}[theorem]{Remark}

\newtheorem{example}[theorem]{Example}
\theoremstyle{remark}
\theoremstyle{remark}

\def\thetheorem{\arabic{theorem}}

\begin{document}

\begin{abstract}
Let $K$ be an abelian extension of the rationals. Let $S(K)$ be the Schur group of $K$ and let $CC(K)$ be the subgroup
of $S(K)$ generated by classes containing cyclic cyclotomic algebras. We characterize when $CC(K)$ has finite index in
$S(K)$ in terms of the relative position of $K$ in the lattice of cyclotomic extensions of the rationals.
\end{abstract}

\maketitle

\section{Introduction}

Throughout this article, $K$ is an abelian extension of the rationals, $Br(K)$ denotes the Brauer group of $K$ and
$S(K)$ the Schur subgroup of $K$. Recall that a {\it cyclotomic algebra} over $K$ is a crossed product $(E/K,\alpha)$,
where $E/K$ is a finite cyclotomic extension and $\alpha$ is a factor set taking values in the group of roots of unity
of $E$. If $(E/K,\alpha)$ is a cyclotomic algebra and the extension $E/K$ is cyclic then we say that $(E/K,\alpha)$ is
a {\it cyclic cyclotomic algebra}. Some properties of cyclic cyclotomic algebras with respect to ring isomorphism were
studied in \cite{HOR1}.

It is well known that every element of $Br(K)$ is represented by a cyclic algebra over $K$ and every element of $S(K)$
is represented by a cyclotomic algebra over $K$ \cite{Yam}. However, in general, not every element of $S(K)$ is
represented by a cyclic cyclotomic algebra. In fact, as we will see in this paper, in general, $S(K)$ is not generated
by classes represented by cyclic cyclotomic algebras.

Let $CC(K)$ denote the subgroup of $S(K)$ generated by classes containing cyclic cyclotomic algebras. In other words
$CC(K)$ is formed by elements of $S(K)$ represented by tensor products of cyclic cyclotomic algebras. The aim of this
article is to study the gap between $S(K)$ and $CC(K)$. More precisely, we give a characterization of when $CC(K)$ has
finite index in $S(K)$ in terms of the relative position of $K$ in the lattice of cyclotomic extensions of the
rationals.

For  every positive integer $n$, let $\zeta_n$ denote a complex primitive $n$-th root of unity. By the Benard-Schacher
Theorem \cite{BS}, $S(K)=\bigoplus_p  S(K)_p$, where $p$ runs over the primes such that $\zeta_p\in K$ and $S(K)_p$
denotes the $p$-primary part of $S(K)$.  Thus $CC(K)$ has finite index in $S(K)$ if and only if $CC(K)_p=CC(K)\cap
S(K)_p$ has finite index in $S(K)_p$ for every prime $p$ with $\zeta_p\in K$. Therefore, we are going to fix a prime
$p$ such that $\zeta_p\in K$ and our main result gives necessary and sufficient conditions for
$[S(K)_p:CC(K)_p]<\infty$, in terms of the Galois group of a certain cyclotomic field $F$ that we are going to
introduce next.

Let $L=\Q(\zeta_m)$ be a minimal cyclotomic field containing $K$, $a$ the minimum positive integer such that
$\zeta_{p^a}\in K$, $s$ the minimum positive integer such that $\zeta_{p^{s}} \in L$ and
 $$b = \left\{ \matriz{{ll}
         s,  & \mbox{ if $p$ is odd or $\zeta_4 \in K$}, \\
         s+v_p([K\cap \Q(\zeta_{p^s}):\Q]) + 2, &
                             \mbox{if } \Gal(K(\zeta_{p^{2a+s}})/K) \mbox{ is not cyclic}, \\
         s+1, & \mbox{otherwise,} } \right.$$
where $v_p:\Q\rightarrow \Z$ denotes the $p$-adic valuation.  Then we let   $\zeta=\zeta_{p^{a+b}}$ and define $F = L(\zeta)$.

The Galois groups of $F$ mentioned above are
    $$\Gamma = \Gal(F/\Q), \quad G = \Gal(F/K), \quad C = \Gal(F/K(\zeta)) \quad \text{and} \quad
    D=\Gal(F/K(\zeta+\zeta\inv)).$$
Notice that $D\ne G$ by the definition of $b$, and if $C\ne D$ then $p^a=2$ and $\rho(\zeta)=\zeta\inv$ for every
$\rho\in D\setminus C$.

We need to fix elements $\rho, \sigma$ of $G$, with $G=\GEN{\rho,\sigma,C}$, such that $D=B\times \GEN{\rho}$ and
$C=B\times \GEN{\rho^2}$ for some subgroup $B$ of $C$ and $G/C=\GEN{\rho C}\times \GEN{\sigma C}$. Furthermore, if
$G/C$ is cyclic (equivalently $C=D$) then we select $\rho=1$ and otherwise $\sigma$ is selected so that
$\sigma(\zeta_4)=\zeta_4$. The existence of such $\rho$ and $\sigma$ in $G$ follows by standard arguments (see
\cite[Lemma 1.4]{Pen1} or \cite[Lemma 3]{HOR2}).

The order of a group element $g$ is denoted by $|g|$. Finally, to every $\psi\in \Gamma$ we associate two non-negative
integers,
    $$d(\psi) = \min\{a, \max\{h \ge 0 :  \psi(\zeta_{p^h})=\zeta_{p^h}\}\} \quad \text{and} \quad
    \nu(\psi) = \max\{0,a-v_p(|\psi G|)\},$$
and a subgroup of $C$:
$$T(\psi) = \{ \eta \in B : \eta^{p^{\nu(\psi)}} \in B^{p^{d(\psi)}} \}.$$

Now we are ready to state our main result.

\begin{theorem}\label{finite}
Let $K$ be an abelian extension of the rationals, $p$ a prime integer and use the above notation.

If $G/C$ is cyclic then the following are equivalent.

\begin{enumerate}
\item $CC(K)_p$ has finite index in $S(K)_p$.

\item For every $\psi\in \Gamma_p$ one has
$\psi^{|\psi G|} \in \displaystyle{\bigcup_{i=0}^{|\sigma C|-1}} \sigma^i T(\psi) $.

\item For every $\psi\in \Gamma_p$ satisfying
$\nu(\psi) < \min\{v_p(\exp B),d(\psi)\}$, one has $\psi^{|\psi G|} \in \displaystyle{\bigcup_{i=0}^{|\sigma C|-1}} \sigma^i T(\psi).$

\end{enumerate}

If $G/C$ is non-cyclic then the following are equivalent:

\begin{enumerate}
\item $CC(K)_2$ has finite index in $S(K)_2$.
\item For every $\psi\in \Gamma_2 \setminus G$, if $d=v_2([K\cap \Q(\zeta):\Q])+2$ then
$$\psi^{|\psi G|} \in \Gal(F/\Q(\zeta_{2^{d+1}})) \cap \left(\bigcup_{i=0}^{|\sigma C|-1} \sigma^i \GEN{ \rho, T(\psi) }\right).$$
\end{enumerate}
\end{theorem}

\def\thetheorem{\thesection.\arabic{theorem}}

Notice that conditions (2) and (3) in Theorem~\ref{finite} can be verified by elementary computations in the Galois
group $\Gamma$.

\section{The subgroup of $S(K)$ generated by cyclic cyclotomic algebras}\label{SecCC}

In this section we provide some information on the structure of $CC(K)_p$. We start by introducing some notation and
recalling some known facts about local information concerning $S(K)$.

The group of roots of unity of a field $E$ is denoted by $W(E)$. If $A$ is a central simple $K$-algebra then $[A]$
denotes the class in the Brauer group of $K$ containing $A$.  By a crossed product algebra we mean an associative
algebra $A = (E/K,\alpha) = \bigoplus_{\pi\in \Gal(E/K)} E u_{\pi}$, where $E/K$ is a finite Galois extension,
$\alpha:\Gal(E/K)\times \Gal(E/K) \rightarrow K^*$ is a map, $\{u_\pi : \pi\in \Gal(E/K)\}$ is an $E$-basis of units of
$A$ and the product in $A$ is determined by the rules: $u_{\pi} u_{\tau} = \alpha_{\pi,\tau} u_{\pi \tau}$ and $u_{\pi
}a=\pi(a)u_{\pi}$. The map $\alpha$ is called the {\it factor set} of the crossed product and the basis $\{u_{\pi} :
\pi\in \Gal(E/K)\}$ is called a {\it crossed section} of the crossed product. Replacing each element $u_{\pi}$ of a
crossed section by $v_{\pi} = \lambda_{\pi} u_{\pi}$, for $\lambda_{\pi}$ a non-zero element of $L$ gives rise to
another crossed section of the crossed product which yields a different factor set. This change of crossed section is
called a {\it diagonal change of basis}.  A {\it cyclotomic algebra} over $K$ is a crossed product algebra
$(E/K,\alpha)$ for which $E/K$ is a finite cyclotomic extension and the factor set $\alpha$ takes values in $W(E)$.
When $\Gal(E/K) = \GEN{\tau}$ is a cyclic group, there exists a diagonal change of basis for which the corresponding
factor set is determined by $u_{\tau}^{|\tau|} = \zeta \in K^*$.  The corresponding crossed product algebra is called a
{\it cyclic algebra}, which we will denote by $(E/K,\zeta)$.

Let $\pr=\{r\in \N : r \mbox{ is prime}\}\cup \{\infty\}$. Given $r\in \pr$, we are going to abuse the notation and
denote by $K_r$ the completion of $K$ at a (any) prime of $K$ dividing $r$. If $E/K$ is a finite Galois extension, one
may assume that the prime of $E$ dividing $r$, used to compute $E_r$, divides the prime of $K$ over $r$, used to
compute $K_r$. We use the classical notation:

\begin{tabular}{rcl}
$e(E/K,r)$ & = & $e(E_r/K_r) = $ ramification index of $E_r/K_r$. \\
$f(E/K,r)$ & = & $f(E_r/K_r) = $ residue degree of $E_r/K_r$. \\
$m_r(A)$ &=& Index of $K_r\otimes_K A$, for a Schur algebra $A$ over $K$.
\end{tabular}

Since $E/K$ is a finite Galois extension and $A$ has uniformly distributed invariants, $e(E/K,r)$, $f(E/K,r)$ and
$m_r(A)$ do not depend on the selection of the prime of $K$ dividing $r$ (see \cite{Serre} and \cite{Ben}).

We also use the following notation, for $\pi\subseteq \pr$ and $r\in \pr$:

\begin{tabular}{rcl}
$S(K,\pi)$ &=& $\{[A] \in S(K) : m_{r}(A)=1, \mbox{ for each } r\in \pr\setminus\pi\}$, \\
$S(K,r)$   &=& $S(K,\{r\})$,\\
$CC(K,\pi)$&=& $CC(K)\cap S(K,\pi)$, \\
$CC(K,r)$  &=& $CC(K)\cap S(K,r)$, \\
$\pr_p$ &=& $\big\{r\in \pr\setminus \{\infty\}: CC(K,\{r,\infty\})_p=CC(K,r)_p\bigoplus CC(K,\infty)_p\big\}$.
\end{tabular}

If $p$ is odd or $\zeta_4\in K$ then $m_{\infty}(A)=1$ for each Schur algebra $A$ and so $\pr_p=\pr\setminus\{\infty\}$.

Finally, if $r$ is odd then we set
$$\nu(r) = \max\{0,a+v_p(e(K(\zeta_r)/K,r))-v_p(|W(K_r)|)\}.$$

\smallskip
The following theorem provides information on the structure of $CC(K)_p$.

\begin{theorem}\label{CCDec}
For every prime $p$ we have
    $$CC(K)_p = \left(\bigoplus_{r\in \pr_p} CC(K,r)_p\right) \bigoplus
    \left(\bigoplus_{r\in \pr\setminus\pr_p} CC(K,\{r,\infty\})_p\right).$$
Let $X_r$ denote the direct summand labelled by $r$ (of either the first or the second kind) in the previous decomposition.
\begin{enumerate}
\item
If $r$ is odd then $X_r$ is cyclic of order $p^{\nu(r)}$ and it is generated by the class of
$(K(\zeta_r)/K,\zeta_{p^a})$.
\item
$X_2$ has order $1$ or $2$ and if it has order $2$ then $p^a=2$ and $X_2$ is generated by the class of
$(K(\zeta_4)/K,-1)$.
\item
If $X_{\infty} \ne 1$, then $p=2$, $K\subseteq \R$, and $X_{\infty}$ has exponent $2$.

\end{enumerate}
\end{theorem}

\begin{proof}
Let $A=(E/K,\xi)$ be a cyclic cyclotomic algebra with $[A]\in S(K)_p$. One may assume without loss of generality that
$\xi\in W(K)_p$. For every subextension $M$ of $E/K$ let $M'$ denote the maximal subextension of $M/K$ of degree a
power of $p$. Since $E/K$ is cyclic, the subextensions of $E/K$ of degree a power of $p$ are linearly ordered. This
implies that $E'=K(\zeta_{r^k})'$ for some prime $r$ and integer $k$. Furthermore, if $\ell=[E:K(\zeta_{r^k})]$,
then $[A^{\otimes \ell}] = [(E/K,\xi^{\ell})] = [(K(\zeta_{r^k})/K,\xi)]$. Since $\ell$ is coprime to $p$, $m_q(A)
= m_q(A^{\otimes \ell}) = m_q(K(\zeta_{r^k})/K,\xi)$, for every $q\in \pr$. If $q\not\in \{r,\infty\}$ then
$K(\zeta_{r^k})/K$ is unramified at $q$ and therefore $m_{q}(A)=1$ \cite[pg. 67, Exercise 16]{Rei}.  Thus  $[A]\in
CC(K,\{r,\infty\})$. This shows that $CC(K)_p = \sum_{r\in \pr\setminus \{\infty\}} CC(K,\{r,\infty\})_p$.

If $\pr\setminus \{\infty\}=\pr$ then this implies that $$CC(K)_p = \bigoplus_{r\in \pr} CC(K,r) = \left(\bigoplus_{r\in
\pr_p} CC(K,r)\right) \bigoplus CC(K,\infty)$$ as wanted. Assume otherwise that $\pr \setminus \{\infty\}\ne \pr_p$. If
$1\ne [A]\in CC(K,\infty)$ then for every $r\in \pr$ and $[B]\in CC(K,\{r,\infty\})\setminus CC(K,\infty)$ one has
$[B]=[A\otimes B]\cdot [A]$ and $[A\otimes B]\in CC(K,r)$. This implies that $CC(K,\{r,\infty\})=CC(K,r)\bigoplus
CC(K,\infty)$, contradicting the hypothesis. Hence $CC(K,\infty)=1$ and then $CC(K)_p = \left(\bigoplus_{r\in \pr_p}
CC(K,r)_p\right) \bigoplus \left(\bigoplus_{r\in \pr\setminus\pr_p} CC(K,\{r,\infty\})_p\right),$ as desired.

(1). Let $r\in \pr$. The map $K_r\otimes_K -: X_r\rightarrow S(K_r)$ is an injective group homomorphism. If
$r$ is odd then $S(K_r)$ is cyclic of order $e(K(\zeta_r)/K,r)$ and it is generated by the cyclic algebra
$(K_r(\zeta_r)/K_r,\zeta_n)$, where $n=|W(K_r)|$ (see e.g. \cite{Yam}). Therefore $X_r$ is cyclic and hence
it is generated by a class containing a cyclic cyclotomic algebra $A$. As above we may assume that
$A=(K(\zeta_{r^k})/K,\zeta_{p^a}^{\ell})$ for some $k,\ell \ge 1$. Since $[A] =
[(K(\zeta_{r^k})/K,\zeta_{p^a})]^{\ell}$, one may assume that $\ell=1$. Then $|X_r| = m_r(A) =
m_r((K(\zeta_{r^k})'/K,\zeta_{p^a})) = m_r((K(\zeta_r)/K,\zeta_{p^a})) = m((K_r(\zeta_r)/K_r,\zeta_{p^a})) =
m((K_r(\zeta_r)/K_r,\zeta_{p^{a+a(r)}})^{\otimes a(r)}) = p^{\nu(r)}$, where $a+a(r)=v_p(n)$. This proves (1).

(2) and (3) follow by similar or standard arguments.
\end{proof}

\begin{remark}\label{CCRemark}\rm
Notice that the proof of Theorem~\ref{CCDec} shows that if $A$ is a cyclic cyclotomic algebra of index a power of $p$ then $[A]\in S(K,\{r,\infty\})$ for some prime $r\in \pr\setminus\{\infty\}$ and if $p$ is odd or $\zeta_4 \in K$ then $[A]\in S(K,r)$.
\end{remark}

By Theorem~\ref{CCDec}, if $r$ is odd then $\nu(r)=\max\{v_p(m_r(A)) : [A] \in CC(K)_p\}$. We can extend the definition
of $\nu(r)$ by setting $\nu(2)=\max\{v_p(m_2(A)) : [A] \in CC(K)_p\}$. Notice that $\nu(2)\le 1$ and $\nu(2)=1$ if and
only if $p^a=2$ and $(K(\zeta_4)/K,-1)$ is non-split.  We will need to compare $\nu(r)$ to
$$ \beta(r) = \max\{v_p(m_r(A)) : [A] \in S(K)_p \}. $$

 A consequence of Theorem~\ref{CCDec} is the following.

\begin{corollary}\label{CC=S}
Let $r\in \pr$. Then
\begin{enumerate}
\item $CC(K)_p=S(K)_p$ if and only if $\nu(r)=\beta(r)$ for each $r\in \pr\setminus\{\infty\}$.
\item $CC(K)_p$ has finite index in $S(K)_p$ if and only if $\nu(r)=\beta(r)$ for all but finitely many primes $r$.
\end{enumerate}
\end{corollary}

\begin{proof}
We prove (2) and let the reader to adapt the proof to show (1).

Assume that $CC(K)_p$ has finite index in $S(K)$ and let $[A_1],\dots,[A_n]$ be a complete set of representatives of
cosets modulo $CC(K)_p$. Then $\pi=\{r\in \pr: m_r(A_i)\ne 1 \text{ for some } i\}$ is finite and $\nu(r)=\beta(r)$ for
every $r\in \pr \setminus \pi$. Conversely, assume that $\nu(r)=\beta(r)$ for every $r\in \pr\setminus \pi$, with $\pi$
a finite subset of $\pr$ containing $\infty$. Then $S(K,\pi)_p$ is finite and we claim that $S(K)_p =
S(K,\pi)_p+CC(K)_p$. Let
$[B]\in S(K)_p$. We prove that $[B]\in S(K,\pi)_p+CC(K)_p$ by induction on $h(B)=\prod_{r\in \pr\setminus\pi}
m_r(B)$. If $h(B)=1$ then $[B]\in S(K,\pi)_p$ and the claim follows. Assume that $h(B) > 1$ and the induction
hypothesis. Then there is a cyclic cyclotomic algebra $A$ and $r\in \pr\setminus \pi$ such that $m_r(B)=m_r(A)>1$.
Since $S(K_r)$ is cyclic, there is a positive integer $\ell$ coprime to $m_r(B)$ such that $(A^{\otimes \ell})
\otimes_K K_r \cong B \otimes_K K_r$ as $K_r$-algebras. Let $C=(A^{\text{op}})^{\otimes \ell} \otimes B$. Since
$A^{\otimes \ell} \in CC(K,\{r,\infty\})_p$, it follows that $h(C)=\frac{h(B)}{m_r(A)}<h(B)$, and hence $[C]\in
S(K,\pi)_p+CC(K)_p$, by the induction hypothesis. Therefore, $[B]=[A]^{\ell}[C]\in S(K,\pi)_p+CC(K)_p$, as required.
\end{proof}

Notice that for $p$ odd Corollary~\ref{CC=S} is a straightforward consequence of the decomposition of $CC(K)_p$ given
in Theorem~\ref{CCDec} and the Janusz Decomposition Theorem \cite{Jan2}.

\section{Examples}

In this section we present several examples comparing $S(K)$ and $CC(K)$ for various fields. We use the standard notation for the generalized quaternion algebra:
$$\quat{a,b}{K}=K[i,j | i^2=a, j^2=b, ji=-ij, a,b\in K^*] \quad \mbox{ and }\quad \HQ(K)=\quat{-1,-1}{K}.$$

\begin{example} $K=\Q$. \rm

It follows from the Hasse--Brauer--Noether--Albert Theorem that $CC(\Q,r)$ is trivial for all primes $r$. The cyclic
cyclotomic algebra $\HQ_{2,\infty} = \HQ(\Q) = (\Q(\zeta_4)/\Q,-1)$ is a rational quaternion algebra which lies in
$CC(\Q,\{2,\infty\})$. When $r$ is odd, the cyclic algebra $\HQ_{r,\infty}=(\Q(\zeta_r)/\Q,-1)$ has real completion
$\R\otimes_{\Q}\HQ_{r,\infty} \simeq M_n(\HQ(\R))$, for $n=\frac{r-1}{2}$, so $m_{\infty}(\HQ_{r,\infty})=2$. The
extension $\Q_r(\zeta_r)/\Q_r$ is unramified at primes other than $r$, so $[\HQ_{r,\infty}]\in CC(\Q,\{r,\infty\})$
(and $m_r(\HQ_{r,\infty})$ must be $2$).  If $r$ and $q$ are distinct finite primes, then
$[\HQ_{r,\infty}][\HQ_{q,\infty}]$ is an element of $CC(\Q,\{r,q\})$, and it follows from Remark~\ref{CCRemark}
that this element cannot be represented by a cyclic cyclotomic algebra. Nevertheless, it is easy to see at this point
that $S(\Q)=CC(\Q)$.

The smallest example of an algebra representing an element in $CC(\Q,\{2,3\})$ is the generalized quaternion algebra
$\quat{-3,2}{\Q}$. The algebra of $2\times 2$ matrices over $\quat{-3,2}{\Q}$ is isomorphic to a simple component of
the rational group algebra of the group of order $48$ that has the following presentation $\langle
x,y,z:x^{12}=y^2=z^2=1, x^y=x^5, x^z=x^7, [y,z]=x^9 \rangle.$
\end{example}

\begin{example} $CC(K,\infty)\ne 1$.

\rm It is also possible that $CC(K,\infty)$ is non-trivial. For example, $\HQ(\Q(\sqrt{2})) = (\Q(\zeta_8)/\Q(\sqrt{2}),-1)$ is homomorphic to a simple component of the rational group algebra of the generalized quaternion group of order $16$. It has real completion isomorphic to $\HQ(\R)$ at both infinite primes of $\Q(\sqrt{2})$, so $m_{\infty}(\HQ(\Q(\sqrt{2}))) = 2$.  If $r$ is an odd prime then $m_r(\HQ(\Q(\sqrt{2})))=1$. Since $\Q_2(\sqrt{2})/\Q_2$ is ramified and the sum of the local invariants at infinite primes is an integer, we deduce that $m_2(\HQ(\Q(\sqrt{2}))) = 1$, so it follows that $[\HQ(\Q(\sqrt{2}))]\in CC(\Q(\sqrt{2}),\infty)$.
\end{example}

\begin{example} Cyclotomic fields.

\rm Suppose $K = \mathbb{Q}(\zeta_m)$ for some positive integer $m > 2$. Assume that either $m$ is odd or $4|m$. The
main theorem of \cite{Jan3} shows that if $p$ is a prime dividing $m$ and $m = p^nm_0$ with $m_0$ coprime to $p$, then
$$S(\mathbb{Q}(\zeta_m))_p = \{ [A \otimes_{\mathbb{Q}(\zeta_{p^n})}\mathbb{Q}(\zeta_m)] :[A] \in S(\mathbb{Q}(\zeta_{p^n}))_p \}.$$
When $p^n > 2$, we know by \cite[Theorem 3]{BS} that $S(\mathbb{Q}(\zeta_{p^n}))_p$ is generated by the Brauer classes of characters of certain metacyclic groups, which, in their most natural crossed product presentation, take the form of cyclic cyclotomic algebras.  Therefore, $S(\mathbb{Q}(\zeta_{p^n}))_p = CC(\mathbb{Q}(\zeta_{p^n}))_p$.  Since it is easy to see that when $K$ is an extension of a field $E$, $\{ [A \otimes_E K] : [A] \in CC(E) \} \subseteq CC(K)$, we can conclude that $S(\mathbb{Q}(\zeta_m)) = CC(\mathbb{Q}(\zeta_m))$ for all positive integers $m$.
\end{example}

Combining Corollary~\ref{CC=S} with the results of \cite{Jan2} one can obtain examples with $S(K)_p\ne CC(K)_p$.

\begin{example}\label{NoDecomp} $CC(K)_p\ne S(K)_p$, $p$ odd.

\rm By Theorem~\ref{CCDec}, if $CC(K)_p=S(K)_p$ then $S(K)_p=\bigoplus_{r\in \pr} S(K,r)_p$. However Proposition 6.2 of
\cite{Jan2} shows that for every odd prime $p$ there are infinitely many abelian extensions $K$ of $\Q$ such that
$S(K)_p\ne \bigoplus_{r\in \pr} S(K,r)_p$. Thus for such fields $K$ one has $S(K)_p\ne CC(K)_p$.
\end{example}

\begin{example}\label{NoDecomp2} $CC(K)_2\ne S(K)_2$ with $\zeta_4\in K$.

\rm Let $q$ be a prime of the form $1+5\cdot 2^9 t$ with $(t,10)=1$. In the last section of \cite{Jan3} one constructs
a subfield $K$ of $\Q(\zeta_{2^9\cdot 5 \cdot q})$ such that $\max\{m_q(A) : [A] \in S(K)_2\}=4$ (in particular
$\zeta_4\in K$), and for every $[A]\in S(K)_2$ with $m_q(A)=4$, one has $m_r(A)\ne 1$, for some prime $r$ not dividing
$10q$. In the notation of Corollary~\ref{CC=S} this means that $v_2(|S(K,q)|)<\beta(q)=4$ (for $p=2)$. Then $S(K)_2 \ne
\bigoplus_{r\in \pr} S(K,r)_2$ and, as in Example~\ref{NoDecomp}, this implies that $CC(K)_2 \ne S(K)_2$.
\end{example}

\begin{example}  $CC(K)_2\ne S(K)_2$ with $\zeta_4\not\in K$.

\rm An example with $S(K)_2\ne CC(K)_2$ and $\zeta_4\not\in K$ can be obtained using Theorem~5 of \cite{Jan1}. This
result gives necessary and sufficient conditions for $S(k)$ to have order $2$ when $k$ is a cyclotomic extension of
$\Q_2$. This is the maximal $2$-local index for a Schur algebra. If $|S(k)|=2$ then $\zeta_4\not\in k$. If, moreover,
$\HQ=(k(\zeta_4)/k,-1)$ is not split then $CC(k)_2=S(k)_2$, because $\HQ$ is a cyclic cyclotomic algebra. However,
there are some fields $k$ for which $|S(k)|=2$ and $\HQ$ is split. In that case $S(k)$ is generated by the class of a
cyclotomic algebra $A$ and we are going to show that $CC(k)\ne S(k)$. Then for any algebraic number field with $K_2=k$
we will also have $CC(K)_2\ne S(K)_2$.

Indeed, if $CC(k)=S(k)$ then $A$ is equivalent to a cyclic cyclotomic algebra $(k(\zeta_m)/k,\zeta)$. One may assume
that $\zeta\in W(k)_2\setminus \{1\}$ and hence $\zeta=-1$, because $\zeta_4\not\in k$. Write $m=2^{v_2(m)}m'$, with
$m'$ odd. Since $k(\zeta_m)/k$ must be ramified, $v_2(m)\ge 2$. If $k(\zeta_{m'})/k$ has even degree then this would
contradict the fact that $k(\zeta_m)/k$ is cyclic. So $k(\zeta_{m'})/k$ has odd degree and therefore
$(k(\zeta_m)/k,-1)$ is equivalent to $(k(\zeta_{2^{v_p(m)}})/k,-1)$ by \cite[(30.10)]{Rei}. Then $A$ is equivalent to
$(k(\zeta_4)/k,-1)$ by \cite[Theorem~1]{Jan1}, yielding a contradiction.
\end{example}

\section{Finiteness of $S(K)_p/CC(K)_p$}

In this section we prove Theorem~\ref{finite}.  The main idea is to compare $\nu(r)$ and $\beta(r)$ for odd primes $r$
not dividing $m$.  We will use the notation introduced in sections 1 and 2 including the Galois groups $\Gamma$, $G$, $C$, $D$, and $B$, the
elements $\rho,\sigma\in G$, and the decompositions $D=\GEN{\rho}\times B$ and $C=\GEN{\rho^2}\times B$.

We also use the following numerical notation for every odd prime $r$ not dividing $m$:

\begin{tabular}{rcl}
$a+a(r)$ &=& $v_p(|W(K_r)|)$, \\
$d(r)$ & = & $\min\{a,v_p(r-1)\}$, \\
$f_r$ & = & $f(K/\Q,r)$, \\
$f(r)$ &=& $v_p(f_r)$, \\
\end{tabular}

\noindent and introduce $\psi_r\in \Gamma$ and $\phi_r\in G$ as follows:
$$\psi_r(\varepsilon)=\varepsilon^r \text{ for every root of unity } \varepsilon \in F, \quad \text{and} \quad \phi_r=\psi_r^{f_r}.$$

The order of $\psi_r$ modulo $G$ is $f_r$, and $\psi_r$ and $\phi_r$ are Frobenius automorphisms at $r$ in $\Gamma$ and
$G$ respectively. By the uniqueness of an unramified extension of a local field of given degree, one has $v_p(|W(K_r)|)
= v_p(|W(\Q_r)|)+f(r)=v_p(e(K(\zeta_r)/K,r))+f(r)$. Thus
    \begin{equation}\label{nu}
    \nu(r) = \max\{0,a-f(r)\}.
    \end{equation}

This gives $\nu(r)$ in terms of the numerical information associated to $r$. Next we quote a result from \cite{HOR2}
which gives the value of $\beta(r)$. This result was obtained by Janusz \cite[Theorem 3]{Jan2} in the case when $p$ is
odd or $\zeta_4 \in K$.  The remaining case was considered by Pendergrass in \cite{Pen1}, but the results there were
based on incorrect calculations involving factor sets (see \cite{HOR2}).

\begin{theorem} \label{localindex}
Let $r$ be an odd prime not dividing $m$ and use the above notation. Let $\phi_r=\psi_r^{f_r} = \rho^{j'} \sigma^j
\eta$, with $\eta \in B$, $0 \le j' < |\rho|$ and $0 \le j < |\sigma C|$.

\begin{enumerate}
\item If $G/C$ is non-cyclic (and hence $p^a=2$) and $j \not\equiv j' \mod 2$, then $\beta(r) = 1$.
\item Otherwise $\beta(r) = \max\{\nu(r),v_p(|\eta B^{p^{d(r)}}|) \}$.
\end{enumerate}
\end{theorem}

We will need the following lemma.

\begin{lemma}\label{alpha}
$\nu(r)$ and $\beta(r)$ depend only on $d(r)$ and the element $\psi_r \in \Gamma$.
\end{lemma}

\begin{proof}
$\nu(r)$ is determined by $f(r)$ (see (\ref{nu})), and  $f(r)$ by $f_r=|\psi_r G|$.  So $\nu(r)$ is determined by
$\psi_r$.  On the other hand, $\psi_r=\rho^{j'} \sigma^j \eta$ for uniquely determined integers $0\le j' <
|\rho|$, $0\le j < |\sigma C|$ and $\eta \in B$.  Therefore, $\psi_r$ determines whether or not $j \equiv j' \mod 2$,
and also the element $\eta$ required in Theorem~\ref{localindex}. So knowing $\psi_r$ and $d(r)$ will allow one to
compute $\beta(r)$.
\end{proof}

We can now give a necessary and sufficient condition, in local terms, for $CC(K)_p$ to have finite index in $S(K)_p$.

\begin{theorem}\label{local1}
$CC(K)_p$ has finite index in $S(K)_p$ if and only if $\nu(r)=\beta(r)$ for all odd primes $r$ not dividing $m$.
\end{theorem}

\begin{proof}
The sufficiency is a consequence of Corollary~\ref{CC=S}.

Suppose that there is an odd prime $r$ not dividing $m$ for which $\nu(r) < \beta(r)$. By Dirichlet's Theorem on primes
in arithmetic progression there are infinitely many primes $r'$ such that $r'\equiv r \mod
\lcm(m,p^{a+b},p^{v_p(r-1)+1})$. For such an $r'$ one has $\psi_{r'} = \psi_r$ and $v_p(r'-1) = v_p(r-1)$. Then
$\beta(r')= \beta(r)>\nu(r)=\nu(r')$ for infinitely many primes $r'$, by Lemma~\ref{alpha}, and hence
$[S(K)_p:CC(K)_p]=\infty$, by Corollary~\ref{CC=S}.
\end{proof}

When $p$ is odd, this result can be interpreted in terms of the local subgroups of $S(K)_p$ and $CC(K)_p$.

\begin{theorem}\label{local}
Let $K$ be a subfield of $\Q(\zeta_n)$ and $p$ an odd prime and $n$ a positive integer. Then the following conditions are equivalent:
\begin{enumerate}
\item $CC(K)_p$ has finite index in $S(K)_p$.
\item $CC(K,r)_p=S(K,r)_p$, for almost all $r\in \pr$.
\item $CC(K,r)_p=S(K,r)_p$, for every prime $r$ not dividing $n$.
\end{enumerate}
\end{theorem}

\begin{proof}
By the Janusz Decomposition Theorem \cite{Jan2}, we have
$$S(K)_p = S(K,\pi)_p \bigoplus \left(\bigoplus_{r \not\in \pi} S(K,r)_p\right),$$
where $\pi$ is the set of prime divisors of $m$, the smallest integer for which $K \subseteq
\mathbb{Q}(\zeta_m)$. This shows that $\beta(r) = v_p(|S(K,r)|_p)$, whenever $r$ is a prime that does not divide $m$
and hence, for such primes $\nu(r)=\beta(r)$ if and only if $CC(K,r)_p=S(K,r)_p$. Now the results follows from
Corollary~\ref{CC=S} and Theorem~\ref{local1}.
\end{proof}

An obvious consequence of Theorem~\ref{local} is the following:

\begin{corollary}
If $K$ is a subfield of $\Q(\zeta_n)$ and $p$ is an odd prime then the order of the group $\bigoplus_{r \in \pr,r \nmid n} S(K,r)_p/CC(K,r)_p$ is either $1$ or infinity.
\end{corollary}

We now proceed with the proof of the main theorem.

\medskip {\sc Proof of Theorem 1.}
For each $\psi\in \Gamma$ we put $h(\psi)=\max\{0\le h \le a+b : \psi(\zeta_{p^h})=\zeta_{p^h}\}$.  Clearly
$d(\psi)=\min\{a,h(\psi)\}$.
By Dirichlet's Theorem on primes in arithmetic progression, for every $\psi\in \Gamma$ there exists an odd prime $r$
not dividing $m$ such that $\psi=\psi_r$.  For such a prime one has $h(\psi)=\min\{a+b,v_p(r-1)\}$.  This prime
$r$ can be selected so that $h(\psi)=v_p(r-1)$, because otherwise we would have  $h(\psi)= a+b < v_p(r-1)$, and we
could replace $r$ by a prime $r'$ satisfying $r' \equiv r \mod m$ and $r' \equiv 1+p^{a+b} \mod p^{a+b+1}$. For
such an $r'$, one has $d(r)=d(r')$, and thus $\nu(r)=\nu(r')$ and $\beta(r)=\beta(r')$ by Lemma~\ref{alpha}.

Let $q=|\sigma C|$.

We now consider the case when $G/C$ is cyclic. Then $D=C=B$ and $\rho=1$. We set $t = v_p(\exp(B))$.  If $t=0$,
then $T(\psi)=B$ for every $\psi\in \Gamma_p$, so that (2) and (3) obviously hold. Furthermore
$|\eta B^{p^{d(r)}}|=1$ and so $\nu(r)=\beta(r)$ for all odd primes $r$ not dividing $m$, by
Theorem~\ref{localindex}. So (1) holds by Theorem~\ref{local1}. So to avoid trivialities we assume that $t>0$.

(1) implies (2). Suppose $K$ does not satisfy condition (2) and let $\psi \in \Gamma_p$ with
$\psi^{|\psi G|} \not\in \bigcup_{j=0}^{q-1} \sigma^i T(\psi)$.
Let $r$ be an odd prime not dividing $m$ for which $\psi=\psi_r$ and $h(\psi)=v_p(r-1)$.
Then $d(r)=d(\psi)$ and $p^{f(r)} = f_r = |\psi G|$, so $\nu(r)=\nu(\psi)$.  The assumption $\psi^{|\psi G|} \not\in \bigcup_{j=0}^{q-1} \sigma^i T(\psi)$ means that when we express $\psi^{|\psi G|}$ as $\sigma^j \eta$ with $0 \le j < q$ and  $\eta \in B$,
the order of $\eta B^{p^{d(\psi)}}$ in $B/B^{p^{d(\psi)}}$ is strictly greater than $p^{\nu(\psi)} = p^{\nu(r)}$.  By
Theorem~\ref{localindex}, we have $\beta(r) > \nu(r)$ for this odd prime $r$ not dividing $m$, and so
Theorem~\ref{local1} implies that (1) fails.

(2) implies (3) is obvious.

(3) implies (1). Assume that (1) fails. Then, by Theorem~\ref{local1}, there exists a prime $r$ not dividing $m$ for
which $\beta(r)>\nu(r)$.  As above, we may select such an $r$ so that $v_p(r-1) \le a+b$.

 Let $\psi=\psi_r$.  Our choice of $r$ implies that $d(\psi)=d(r)$.  We claim that one can assume $\psi \in
\Gamma_p$. If $\psi \not\in \Gamma_p$, then let $\ell$ be the least positive integer such that $\psi^{\ell}$ lies in
$\Gamma_p$.  Let $r'$ be a prime integer such that $r'\equiv r^{\ell} \mod \lcm(m,p^{a+b})$.  Since $\ell$ is coprime to
$p$, we have $v_p(r'-1) = v_p(r^{\ell}-1) = v_p(r-1)$ and therefore $d(\psi^{\ell})=d(r')=d(r)=d(\psi)$.  Since
$\psi_{r'}=\psi_r^{\ell}$ and $\ell$ is coprime to $p$, we also have $f(r')=f(r)=f(\psi)$. It follows from Lemma~\ref{alpha}
that $\beta(r)=\beta(r')$ and $\nu(r)=\nu(r')$. So by replacing $r$ by $r'$ if necessary, one may assume that
$\psi\in\Gamma_p$ and $d(\psi)=d(r)$.

For this prime $r$ and element $\psi =\psi_r \in \Gamma_p$, the assumption $\beta(r) > \nu(r)$ and Theorem~\ref{localindex} imply that, when we write $\phi_r=\psi^{f_r} = \sigma^j \eta$, with $0 \le j < q$ and $\eta \in B$,
the order of $\eta B^{p^{d(r)}}$ in $B/B^{p^{d(r)}}$ is precisely $p^{\beta(r)}$. Then $\eta^{p^{\nu(r)}}\not\in B^{p^{d(r)}}$, equivalently $\eta\not\in T(\psi)$ and hence $\psi^{|\psi G|} \not\in \bigcup_{j=0}^{q-1} \sigma^i T(\psi)$.

Since the exponent of $B/B^{p^{d(r)}}$ is precisely $p^k$, where $k = \min\{t,d(r)\}$, this can only be possible if
$\nu(\psi)=\nu(r) < k = \min\{t,d(\psi)\}$. This shows that if condition (1) fails, then condition (3) also
fails.  This completes the proof in the case that $G/C$ is cyclic.

Now suppose $G/C$ is non-cyclic.  In particular, $p^a=2$ and $\sigma(\zeta_4)=\zeta_4$. Let $d=v_2([K\cap \Q(\zeta):\Q])+2$ and let $c$ be an integer such that $\sigma(\zeta)=\zeta^c$. Then $v_p(c-1)=d$ and $d(\psi)=1$ for all $\psi \in \Gamma_2$.

(1) implies (2). Suppose (2) fails. Then there exists a $\psi \in \Gamma_2 \setminus G$ such that either
$\psi^{|\psi G|} \not\in \Gal(F/\Q(\zeta_{2^{d+1}}))$ or
$\psi^{|\psi G|} \not\in \bigcup_{i=0}^{q-1} \sigma^i\GEN{ \rho, T(\psi) }$.

As above, there exists an odd prime $r$ not dividing $m$ such that $\psi = \psi_r$, $|\psi G| = f_r$, and
$\nu(\psi)=\nu(r)$. Since $\psi \not\in G$ we have $f(r) > 0$ and so $\nu(r) = 0$, by (\ref{nu}). Also from
$f(r)>0$ one deduces that $\phi_r=\psi^{f_r}$ fixes $\zeta_4$ and so when we write $\phi_r=\rho^{j'}\sigma^j \eta$ with $0\le j'<
|\rho|$, $0\le j \le q$, $\eta\in B$, we have that $j'$ is even.

If $\phi_r \not\in \Gal(F/\Q(\zeta_{2^{d+1}}))$, then $j$ is odd, and we are in the case of Theorem~\ref{localindex},
part (1), with $\nu(r)=0$ and $\beta(r)=1$. Otherwise $j$ is even and $\psi^{|\psi G|} \not\in \bigcup_{i=0}^{q-1} \sigma^i\GEN{ \rho, T(\psi) }$.
Then $\eta\not\in T(\psi)$, or equivalently, $\nu(r)<v_2(|\eta B^2|)$ (observe that $d(r)=1$). By Theorem~\ref{localindex}, we have
$\beta(r)=v_2(|\eta B^2|)>\nu(r)$.
Therefore, in all cases in which (2) fails, we have $\nu(r)<\beta(r)$.  So (1) fails by
Theorem~\ref{local1}.

(2) implies (1). Suppose (1) fails.  By Theorem~\ref{local1}, there exists an odd prime $r$ not dividing $m$ such that
$0 = \nu(r) < \beta(r) = 1$.  Since $\nu(r) = 0$, we must have $f(r) > 0$,  so $\psi = \psi_r \not\in G$.  As
above, we may adjust $\psi_r$ by an odd power and make a different choice of $r$ without changing $\nu(r)$ or
$\beta(r)$ in order to arrange that  $\psi \in \Gamma_2$.  Write $\phi_r=\psi^{f_r} = \psi^{|\psi G|} =
\rho^{j'}\sigma^j \eta$, with $0\le j' < |\rho|$, $0 \le j < q$ and $\eta \in B$. As above, $j'$ is even
because $f(r)>0$. If $j$ is odd, then $\psi^{|\psi G|}\not\in \Gal(F/\Q(\zeta_{2^{d+1}}))$ and so (2) fails.
Suppose now that $j$ is even, so we have $\psi^{|\psi G|} \in \Gal(F/\Q(\zeta_{2^{d+1}}))$. Then the fact that
$\beta(r)=1$ implies by Theorem~\ref{localindex}, part (2), that $|\eta B^{2^{d(r)}}|=2$. Since $d(r)\le a=1$,
we have $d(\psi)=d(r)=2$ and so $\eta\not\in B^2$ and $\eta\not\in T(\psi)$. Then $\psi^{|\psi G|} \not\in  \prod_{i=0}^{q-1} \sigma^i \GEN{\rho,T(\psi)}$ and so (2) fails. \hfill $\square$

\medskip
Some obvious consequences of Theorem~\ref{finite} are the following.

\begin{corollary}\label{SuffFinite2}
If $\psi^{|\psi G|} \not\in \GEN{\sigma,\rho,T(\psi)}$, for some $\psi\in \Gamma_p$, then $CC(K)_p$ does not have finite index in $S(K)_p$.
\end{corollary}

\begin{corollary}\label{SuffFinite}
If $G/C$ is cyclic and $\nu(\psi) \ge \min\{v_p(\exp B),d(\psi)\}$ for all $\psi\in \Gamma_p$, then $CC(K)_p$ has finite index in $S(K)_p$.
\end{corollary}

\begin{corollary}\label{t=1b=0}
If $G/C$ is cyclic and $v_p(\exp B)+v_p(\exp(\Gal(K/\Q)))\le a$ then $CC(K)_p$ has finite index in $S(K)_p$.
\end{corollary}

\begin{proof}
If $\psi\in \Gamma_p$ then $v_p(|\psi G|) \le v_p(\exp(\Gal(K/\Q)))\le a-v_p(\exp B)$, by assumption. Therefore $\nu(\psi)=\max\{0,a-v_p(|\psi G|)\} \ge v_p(\exp B)$ and Corollary~\ref{SuffFinite} applies.
\end{proof}

\begin{example} A simple example with $[S(K)_p:CC(K)_p]=\infty$. \rm

Let $p$ and $q$ be odd primes with $v_p(q-1)=2$. Let $K$ be the subextension of $L=\Q(\zeta_{pq})/\Q(\zeta_p)$ with
index $p$ in $\Q(\zeta_{pq})$. Then $F=\Q(\zeta_{p^2q})$, $G \cong \langle \theta \rangle \times C$ is elementary
abelian of order $p^2$, and $\Gamma_p$ has an element $\psi$ such that $\psi^p$ generates $C$. Then $a=v_p(|\psi G|)=1$
and so $\nu(\psi)=0$ and $d(\psi)=1$. Therefore, $T(\psi)=1$ and hence $\GEN{\sigma,T(\psi)}=\GEN{\sigma}$. However
$\GEN{\sigma}\cap C = 1$ and hence $\psi^{|\psi G|}=\psi^p \not \in \GEN{\sigma,T(\psi)}$.  So it follows from
Corollary~\ref{SuffFinite2} that $CC(K)_p$ has infinite index in $S(K)_p$.
\end{example}

The reader may check using Theorem~\ref{finite} that $[S(K)_p:CC(K)_p]=\infty$ for the fields $K$ constructed by Janusz that were mentioned in Example~\ref{NoDecomp}. The same holds for the field of Example~\ref{NoDecomp2}. This can be verified using the arguments in the proofs of Lemmas 4.2 and 6.4 and Proposition 6.5 in \cite{Jan2} where it is proved that $0=v_p(|S(K,q)|)<\beta(q)$ for all the primes $q$ such that $q\equiv 1 \mod 16$ and $r$ is not a square modulo $q$.

In all the examples shown so far, the index of $CC(K)_p$ in $S(K)_p$ is
either 1 or infinity.  This, together with Corollary 4.5, may lead one to
believe that the quotient group $S(K)_p/CC(K)_p$ is either trivial or
infinite for every field $K$ and every prime $p$. By Corollary 2.3 and
Theorem 4.3, $S(K)_p/CC(K)_p$ is both finite and non-trivial if and
only if $\nu(r)=\beta(r)$ for every odd prime not dividing $m$ and
$\nu(r)\ne \beta(r)$ for $r$ either 2 or an odd prime dividing $m$. In
the following example we show that for every odd prime $p$ there exists
a field $K$ satisfying these conditions.

\begin{example} An example with $CC(K)_p\ne S(K)_p$ and $[S(K)_p:CC(K)_p]<\infty$.
\rm

Let $p$ be an arbitrary odd prime and let $q$ and $r$ be primes for which $v_p(q-1) = v_p(r-1) = 2$, $v_q(r^p-1)=0$, and $v_q(r^{p^2}-1)=1$.  The existence of such primes $q$ and $r$ for each odd prime $p$ is a consequence of Dirichlet's Theorem on primes in arithmetic progression. Indeed, given $p$ and $q$ primes with $v_p(q-1)=2$, there is an integer $k$, coprime to $q$ such that the order of $k$ modulo $q^2$ is $p^2$. Choose a prime $r$ for which $r\equiv k+q \mod q^2$ and $r\equiv 1 + p^2 \mod p^3$. Then $p$, $q$ and $r$ satisfy the given conditions.

Let $K$ be the compositum of $K'$ and $K''$, the unique subextensions of index $p$ in
$\Q(\zeta_{p^2q})/\Q(\zeta_{p^2})$ and $\Q(\zeta_{p^2r})/\Q(\zeta_{p^2})$ respectively. Then $m=p^2rq$, $a=2$ and $L=\Q(\zeta_m)=K(\zeta_q)\otimes_K K(\zeta_r)$. Therefore, $F = \Q(\zeta_{p^4qr})$, and $G = \Gal(F/K(\zeta_{qr})) \times \Gal(F/K(\zeta_{p^4q})) \times \Gal(F/K(\zeta_{p^4r}))$.  We may choose $\sigma$ so that $\langle \sigma \rangle = \Gal(F/K(\zeta_{qr})) \cong G/C$ has order $p^2$.  The inertia subgroup of $r$ in $G$ is $\Gal(F/K(\zeta_{p^4q}))$, which is generated by an element $\theta$ of order $p$.  Note that $B=C$ and $v_p(\exp(\Gal(K/\Q)))=v_p(\exp B)=1<a=2$.  Hence $K$ satisfies the conditions
of Corollary~\ref{t=1b=0} and so $CC(K)_p$ has finite index in $S(K)_p$.

Since $K=K'\otimes_{\Q(\zeta_{p^2})} K''$ and $K''/\Q(\zeta_{p^2})$ is totally ramified at $r$, we have that $K'_r$ is the maximal unramified extension of $K_r/\Q_r$. It follows from $v_q(r^{p^2}-1) = 1$ and $v_q(r^p-1) = 0$ that $[\Q_r(\zeta_q):\Q_r] = p^2$, and so $[K_r':\Q_r] = p = f(K/\Q,r)$. Therefore $v_p(|W(K_r)|) = v_p(|W(\Q_r)|) + f(r) = v_p(r-1) + 1 = 3$, and so we have $\nu(r) = \max\{0,a+v_p(|\theta|)-v_p(|W(K_r)|)\} = 0$.

Let $\psi_r$ be the Frobenius automorphism of $r$ in $\Gal(F/\Q)$.  Then $\psi_r^p = \sigma^p \eta$, where $\eta \in B$
generates $\Gal(F/K(\zeta_{p^4r}))$.  Since $\langle \theta \rangle \cap \langle \eta \rangle = 1$, there exists a
skew pairing $\Psi:B \times B \rightarrow W(K)_p$ such that $\Psi(\theta,\eta)$ has order $p$.  By \cite[Theorem
13]{HOR2}, it follows that $\beta(r) \ge 1$, and so $S(K)_p \ne CC(K)_p$.
\end{example}

The last example shows that when $G/C$ is noncyclic, it is possible for $CC(K)_2$ to have infinite index in $S(K)_2$ even when $t = v_2(\exp B) = 0$.  It also is a counterexample to \cite[Theorem 2.2]{Pen1}.

\begin{example} \rm An example with $[S(K)_2:CC(K)_2] = \infty$ and $C=1$.

Let $q$ be an odd prime greater than $5$.  Let $K = \Q(\zeta_q, \sqrt{2})$, and consider $[S(K)_2:CC(K)_2]$.  It is easy to see that $p^a=2$ and $b=2$.  Since $\Gal(K(\zeta_{2^4})/K)$ is noncyclic, we compute $s = a+b+v_2([\Q(\sqrt{2}):\Q] + 2 = 6$, and so $F = \Q(\zeta_{64q})$.  Since $\Q(\zeta_q) \subset K$, we have $C = \Gal(F/K(\zeta_{64})) = 1$.  For our generators of $\Gal(F/K)$, we may choose $\rho, \sigma$ such that $\rho(\zeta_q) = \zeta_q$, $\rho(\zeta_{64}) = \zeta_{64}^{-1}$, $\sigma(\zeta_q) = \zeta_q$, and $\sigma(\zeta_{64}) = \zeta_{64}^9$.  Let $r$ be any prime for which $r^2 \equiv 1 \mod q$ and $r \equiv 5 \mod 2^6$.  Then $\psi_r \not\in G$, but $5^2 \equiv 9^3 \mod 64$ implies that $\psi_r^2 = \sigma^3$.  This means that we are in the case of Theorem~\ref{localindex} where $\nu(r)=0$ and $j$ is odd, so $\beta(r)=1$.  So $[S(K)_2:CC(K)_2]$ is infinite.
\end{example}

\end{document}